\documentclass[11pt]{article}
\usepackage{mathrsfs}
\usepackage{latexsym,lineno}
\usepackage{epsfig}
\usepackage{color}
\usepackage{amsmath}\usepackage{fleqn}\usepackage{verbatim}\usepackage{epsf}
\usepackage{amsthm}\usepackage{graphicx, float}\usepackage{graphicx}
\usepackage{amsfonts}\usepackage{amssymb}\usepackage{graphpap}
\usepackage{epic}\usepackage{curves}

\newtheorem{theorem}{Theorem}[section]

\newtheorem{corollary}{Corollary}[section]
\newtheorem{lemma}{Lemma}[section]

\title{\bf
Coloring clique-hypergraph of $K_5$-minor-free graphs\footnote{{\it E-mail addresses}: efshan@shu.edu.cn (E. Shan), lykang@shu.edu.cn (L. Kang)}}
\author{Erfang Shan$^{a,b}$\thanks{\em Corresponding author.}, \, Yuxiao Sun$^{b}$, \, Liying Kang$^{b}$ \\
{\small$^a$School of Management, Shanghai University, Shanghai
200444, China}\\
{\small $^b$Department of Mathematics, Shanghai University, Shanghai
200444, China}}

\date{}
\begin{document}
\maketitle
\begin{abstract}

  A clique-coloring of a graph $G$ is a coloring of the vertices of $G$ so
that no maximal clique of size at least two is monochromatic.
The clique-hypergraph, $\mathcal{H}(G)$, of a graph $G$ has $V(G)$ as its set of vertices and the maximal cliques of $G$ as its hyperedges. A (vertex) coloring of $\mathcal{H}(G)$ is a clique-coloring of $G$.
The clique-chromatic number of $G$ is the least number of colors for which  $G$ admits
a clique-coloring.
  Every planar graph has been proved to be 3-clique-colorable  (Electr. J. Combin. 6 (1999), \#R26). Recently, we showed that every claw-free planar graph, different from an odd
cycle,  is $2$-clique-colorable (European J. Combin. 36 (2014) 367-376).
 In this paper  we  generalize these results to \{claw, $K_5$-minor\}-free graphs.


\end{abstract}
\bigskip

\section{Introduction}
A {\em hypergraph} $\mathcal{H}$ is a pair $(V, \mathcal{E})$ where $V$ is a finite set of
vertices and $\mathcal{E}$ is a family
of non-empty subsets of $V$ called {\em hyperedges}. A $k$-{\em coloring} of $\mathcal{H}$ is a function
$\phi: V\rightarrow \{1, 2, \ldots, k\}$ such that for each $S\in \mathcal{E}$, with $|S|\ge 2$, there exist
$u, v\in S$ with $\phi(u)\neq \phi(v)$, that is, there is no monochromatic hyperedge of size at least two. If such a function exists we say that $\mathcal{H}$ is $k$-{\em colorable}.
The {\em chromatic number} $\chi(\mathcal{H})$ of $\mathcal{H}$ is the smallest $k$ for which $\mathcal{H}$ admits a $k$-coloring. In other words,
a $k$-coloration of $\mathcal{H}$ is a partition $\mathcal{P}$ of $V$ into at most $k$ parts such that no hyperedge of
cardinality at least 2 is contained in some $P\in \mathcal{P}$.

Here we consider hypergraphs arising from graphs: for an undirected simple graph $G$, we call {\em clique-hypergraph} of $G$ (or {\em hypergraph of maximal cliques} of
$G$) the hypergraph $\mathcal{H}(G)=(V(G), \mathcal{E})$ which has the same vertices as $G$ and whose
{\em hyperedges} are the {\em maximal cliques} of $G$ (a {\em clique} is a complete induced subgraph
of $G$, and it is {\em maximal} if it is not properly contained in any other clique).
A $k$-coloring of $\mathcal{H}(G)$ is also called a $k$-{\em clique-coloring} of $G$, and the chromatic
number $\chi(\mathcal{H}(G))$ of $\mathcal{H}(G)$ is called the {\em clique-chromatic number} of $G$, denoted by $\chi_{C}(G)$. If $\mathcal{H}(G)$ is $k$-colorable we say that $G$ is $k$-{\em clique-colorable}.

Note that what we call $k$-clique-coloration here is also called {\em weak k-coloring} by Andreae, Schughart
and Tuza in \cite{as,bz}
 or {\em strong $k$-division} by
Ho\'ang and McDiarmid in \cite{hm}. Clearly, any (vertex) $k$-coloring of $G$ is a $k$-clique-coloring of $G$, so $\chi_{C}(G)\leq \chi(G)$. On the other hand, note that if $G$ is triangle-free (contains no a clique on three vertices), then $\mathcal{H}(G)=G$, which implies
$\chi_C(G))=\chi(G)$. Since the chromatic number of triangle-free graphs is known to be
unbounded \cite{my}, we get that the same is true
for the clique-chromatic number.

The clique-hypergraph coloring problem was posed by Duffus et al. \cite{duff}. In general, clique-coloring can be a very different problem from ordinary vertex coloring \cite{bg}. Clique-coloring is harder than ordinary vertex coloring: it is coNP-complete even to check whether
a 2-clique-coloring is valid \cite{bg}. The complexity of 2-clique-colorability is investigated in \cite{kt},
where they show that it is NP-hard to decide whether a perfect graph is 2-clique-colorable.
However, it is not clear whether this problem belongs to NP.  Recently, Marx  \cite{ma} prove that  it is
$\sum_2^p$-complete to check whether a graph
is 2-clique-colorable.
 On the other hand,  Bacs\'o et al. \cite{bg} proved that almost all perfect graphs are 3-clique-colorable.
 A necessary and sufficient condition for $\chi_{C}\leq k$  on line graphs was given \cite{as}. Recently,
 Campos et al. \cite{cdm} showed that powers of cycles is $2$-clique-colorable, except for odd cycles of size
at least five, that need three colours, and showed that odd-seq circulant graphs are $4$-clique-colorable.
 Many papers focus on finding the classes of graphs with $\chi_{C}=2$. Claw-free perfect graphs and claw-free graphs without an odd  hole are  $2$-clique-colorable \cite{bg}.   Claw-free graphs  of maximum degree at most four, other than an odd  cycle, are  $2$-clique-colorable \cite{bz}.
Many subclasses of odd-hole-free graphs have been studied and showed
to be $2$-clique-colorable \cite{de, de1, duff}.  Other works
considering the clique-hypergraph coloring problem in classes of graphs can be found
in the literature \cite{ghm,hm,jt,lon}.

For planar graphs, Mohar and $\check{S}$krekovski \cite{ms} have shown that
every planar graph is 3-clique-colorable, and Kratochv\'il and Tuza \cite{kt} proposed a
polynomial-time algorithm to decide if a planar graph is 2-clique-colorable (the set of
cliques is given in the input).

Mohar and $\check{S}$krekovski \cite{ms} proved the following theorem.

\begin{theorem}[Mohar and $\check{S}$krekovski \cite{ms}]
 Every planar graph is 3-clique-colorable. \label{thm1.1}
\end{theorem}

Recently, we proved the following result in \cite{slk}.

\begin{theorem}[Shan, Liang and Kang \cite{slk}]
Every claw-free planar graph, different from an odd cycle, is 2-clique-colorable.  \label{thm1.2}
\end{theorem}

The purpose of this paper is to  generalize the above results to $K_5$-minor-free graphs. Section 2 gives
some notation and terminology. In Section 3, we first show that
every edge-maximal $K_5$-minor-free graph
is 3-clique-colorable and every edge-maximal $K_4$-minor-free graph is 2-clique-colorable.
Secondly, we prove that every \{claw, $K_5$-minor\}-free graph $G$, different from an odd cycle, is $2$-clique-colorable and a 2-clique-coloring can be found in polynomial time.

\section{Preliminaries}
Let $G$ be an undirected simple graph with {\em vertex set} $V(G)$ and {\em edge set} $E(G)$.
If $H$ is a subgraph of $G$, then the vertex set of $H$ is denoted by $V(H)$. For $v\in V(G)$,
the {\em open neighborhood} $N(v)$ of  $v$ is $\{u\colon\, uv\in E(G)\}$, and the {\em closed
neighborhood} $N[v]$ of $v$ is  $N(v)\cup\{v\}$. The {\it degree} of the vertex $v$, written $d_G(v)$ or simply $d(v)$,
is the number of edges incident to $v$, that is, $d_G(v)=|N(v)|$.
The maximum and minimum degrees of $G$ are denoted by $\varDelta(G)$
and $\delta(G)$, respectively.
For a subset $S\subseteq V(G)$, the subgraph induced by $S$ is denoted by $G|S$.
As usual, $K_{m,n}$ denotes a complete bipartite graph with  classes of
cardinality $m$ and $n$; $K_n$  is the  complete graph on $n$ vertices, and  $C_n$
is the cycle on $n$ vertices.
The graph $K_{1,3}$ is also called a
{\em claw}, and $K_3$ a {\em triangle}. The graph
$K_4-e$ (obtained from $K_4$ by deleting one edge) is called a {\em diamond}.
A graph $G$ is {\it claw-free} if it does not
contain $K_{1,3}$ as an induced subgraph.
A graph $H$ is a {\em minor} of a graph $G$ if $H$ can be obtained from $G$ by
deleting edges, deleting vertices, and contracting edges. A graph $G$ is $H$-{\em minor-free}, if $G$ has no minor which is isomorphic to
$H$.
The family of $K_5$-minor-free graphs is a generalization of the planar graphs.
For a family  $\{F_1, \ldots, F_k\}$ of graphs, we say that $G$
is $\{F_1, \ldots, F_k\}$-free if it is $F_i$-{\em free} for all $i$.

For an  integer of $k$, a clique of size $k$ of a graph $G$ is called a $k$-{\em clique} of $G$.
The largest such $k$ is  the {\em clique number} of $G$, denoted $\omega(G)$.
A subset $I$ of vertices of  $G$ is called an {\em independent set}
of $G$ if no two vertices of $I$ are adjacent in $G$. The maximum cardinality of an independent set of $G$ is
the {\em independence number} $\alpha(G)$ of $G$.
A set $D\subseteq V(G)$  is called a {\em clique-transversal set} of
$G$ if $D$ meets all cliques of $G$, i.e., $D\cap V(C) \neq
\emptyset$ for every clique $C$ of $G$. The {\em clique-transversal
number}, denoted by $\tau_C(G)$, is  the  cardinality of a
minimum clique-transversal set of $G$.
The notion of clique-transversal set in graphs can be regarded as a
special case of the transversal set in hypergraph theory. Erd\H{o}s et al. \cite{egt} have proved that the problem of finding a minimum clique-transversal set for a graph is NP-hard. It is therefore of
interest to determine bounds on the clique-transversal number of a graph. In \cite{egt} Erd\H{o}s et al. proposed to find sharp estimates on the clique-transversal number $\tau_C$ for particular classes of
graphs (planar graphs, perfect graphs, etc.).

We call $G$ a {\em plane triangulation} if every face of a planar $G$ (including the outer
triangulation face) is bounded by a triangle. Let $G$ be a planar graph and $C$
 a cycle of $G$.   The {\em interior} Int($C$) of $C$  denotes the subgraph of $G$ consisting of $C$ and
all vertices and edges in the disk bounded by $C$. Similarly, Ext$(C)\subseteq G$ is the
{\em exterior} of $C$. Obviously, Int$(C)\cap $Ext$(C)=C$.

\section{$K_5$-Minor-free graphs}

In this section, we first show that every edge-maximal $K_4$-minor-free graph is 2-clique-colorable
and every edge-maximal $K_5$-minor-free graph
is 3-clique-colorable. Secondly, we show that every \{claw, $K_5$-minor\}-free graph, different from an odd cycle, is $2$-clique-colorable. As an immediate corollary, we prove
that every \{claw, $K_5$-minor\}-free graph, different from an odd cycle, has the clique-transversal number
 bounded above by half of its order.

\begin{lemma}{\rm (\cite{di})}\label{lem3.1}
A graph with at least 3 vertices is edge-maximal
without a $K_4$-minor if and only if it can be constructed recursively from
triangles by pasting along $K_2$'s.
\end{lemma}
\begin{theorem}\label{thm3.1}
Every edge-maximal $K_4$-minor-free graph is 2-clique-colorable.
\end{theorem}
\noindent{\bf Proof.} Let $G$ be an edge-maximal $K_4$-minor-free graph. For $|V(G)|=2$, the assertion is trivial.
So we may assume that $|V(G)|\ge 3$. Suppose that $A$ is a subgraph isomorphic to $K_2$ in $G$
and $\phi$ is a (not necessarily proper) coloring  of $A$. We show by induction on $|V(G)|$ that $\phi$  can be extended to a 2-clique-coloring of $G$.
For $|V(G)|=3$, since $G$ is a triangle, the assertion is obvoius.
For $|V(G)|\ge 4$, by Lemma \ref{lem3.1}, we have $G=G_1\cup G_2$ such that $G_1\cap G_2\cong K_2$ where $G_1, G_2$ are proper subgraphs of $G$. Clearly $A$ is a subgraph of  $G_1$ or $G_2$.
Without loss of generality, let $A$ be a subgraph of  $G_1$.   By the
induction hypothesis applied to $G_1$, $\phi$  can be extended to a 2-clique-coloring of $G_1$.
This coloring induces a (not necessarily proper) coloring of $G_1\cap G_2$, and by the
induction hypothesis applied to $G_2$, the coloring of $G_1\cap G_2$
further can be extended to a 2-clique-coloring of  $G_2$. The union of these 2-clique-colorings
of $G_1$ and $G_2$ forms a  2-clique-coloring of $G$. The assertion follows.
~$\Box$

\noindent{\bf Remark 1.} The condition ``edge-maximal" in Theorem \ref{thm3.1} is the best possible, because the graph exhibited in Figure 5 is a $K_4$-minor-free graph and we see that it is not 2-clique-colorable.

Next we shall show that every edge-maximal $K_5$-minor-free graph is 3-clique-colorable.
For this purpose, we need the following lemmas.

\begin{lemma}{\rm ($\check{S}$krekovski \cite{skr})}\label{lem3.2}
Every $K_5$-minor-free graph is 5-choosable.
\end{lemma}

\begin{lemma}\label{lem3.3}
Let $G$ be a $K_5$-minor-free graph with at least one edge such that each
edge of $G$ is contained in some triangle of $G$. Then $G$ has a 3-clique-coloring such that
 no triangle of $G$ is monochromatic.
\end{lemma}
\noindent{\bf Proof.}\, By Lemma \ref{lem3.2}, there is a 5-coloring $\phi$ of $G$. For
$i= 1, \ldots, 5$, let $V_i\subseteq V(G)$ be the set of vertices colored $i$. Now, let $\phi(v)=1$
if $v\in V_1\cap V_2$, let $\phi(v)=2$ if $v\in V_3\cup V_4$, and let $\phi(v)=3$ if $v\in V_5$. Since every maximal clique $K$ in $G$ contains at least 3 vertices, $K$ uses at least 3 colors in the 5-coloring of
$G$, and hence $\phi$ uses at least both colors on $K$. Therefore, $\phi$ is a 3-clique-coloring of $G$ and no triangle of $G$ is monochromatic.
~$\Box$

\begin{lemma}{\rm (Mohar and $\check{S}$krekovski \cite{ms})}\label{lem3.4}
Let $G$ be a connected planar graph whose outer cycle $C$ is a
triangle. Let $\phi: V(C)\rightarrow \{1, 2, 3\}$ be a coloring of $\mathcal{H}(C)$. Then $\phi$ can be
extended to a  3-clique-coloring of $G$ and no triangle of $G$ is monochromatic.
\end{lemma}

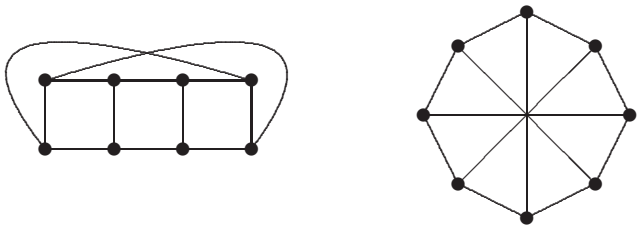
\begin{figure}[ht]\label{map1}
\setlength{\unitlength}{.045in}
\begin{picture}(55,30)
\put(42,-10){\begin{picture}(60,30)

\multiput(-4,28)(8,0){4}{\circle*{1.5}} \multiput(-4,20)(8,0){4}{\circle*{1.5}}
\multiput(-4,28)(8,0){3}{\line(1,0){8}}\multiput(-4,20)(8,0){3}{\line(1,0){8}}
\multiput(-4,28)(8,0){4}{\line(0,-1){8}}
\qbezier(-4,28)(35,40)(20,20)
\qbezier(-4,20)(-20,40)(20,28)

\multiput(40,24)(24,0){2}{\circle*{1.5}}
\multiput(52,12)(0,24){2}{\circle*{1.5}}

\multiput(44,32)(16,0){2}{\circle*{1.5}}
\multiput(44,16)(16,0){2}{\circle*{1.5}}
\put(52,12){\line(0,1){23.5}}\put(40,24){\line(1,0){23.5}}
\put(44,32){\line(1,-1){16}}
\put(44,16){\line(1,1){16}}
\qbezier(40,24)(42,28)(44,32)\qbezier(40,24)(42,20)(44,16)
\qbezier(64,24)(62,28)(60,32)\qbezier(64,24)(62,20)(60,16)
\qbezier(52,36)(48,34)(44,32)\qbezier(52,36)(56,34)(60,32)
\qbezier(52,12)(48,14)(44,16)\qbezier(52,12)(56,14)(60,16)

\end{picture}}
\end{picture}
\caption{Two different representations of the Wagner graph.}
\end{figure}

The {\em Wagner graph}, denoted by $V_8$, is a graph constructed from an 8-cycle (we call it the outer cycle) by connecting the antipodal vertices (these edges will be called the {\em diagonal edges}). The Wagner graph is depicted in Figure 1. Note that the Wagner graph is triangle-free and 3-colorable (because it is cubic).

The following easy lemma about the Wagner graph is obtained by Naserasr et al. in \cite{nns}.
\begin{lemma}{\rm (\cite{nns})}\label{lem3.5}
If $e$ is an edge of the Wagner graph $V_8$, then $V_8-e$ admits a 3-coloring such that the end vertices of $e$ receive the same color.
\end{lemma}

\begin{theorem}{\rm (Wagner \cite{wa})}\label{thm3.2}
If $G$ is an edge-maximal $K_5$-minor-free graph with at least 4 vertices, then $G$ can be constructed recursively, by
pasting along $K_2$'s and $K_3$'s, from plane triangulations and copies of the Wagner graph.
\end{theorem}

In order to make our arguments easier to follow we use the following notations in \cite{nns}: Let $\mathcal{T}=T_1, T_2, \ldots, T_r$ be a sequence of graphs where each $T_i$ is either a plane triangulation or a copy of the Wagner graph $V_8$. By $\mathcal{T}$, we construct another sequence $\mathcal{G}=G_1, G_2, \ldots, G_r$ of graphs as follows: $G_1=T_1$, $G_i$ is obtained from $G_{i-1}$ and $T_i$ by pasting $T_i$ to $G_{i-1}$ along a $K_2$ or a $K_3$.

Given an edge-maximal $K_5$-minor-free graph $G$, the sequence $\mathcal{T}$ is said to be a {\em Wagner sequence} of the graph $G$, if $G=G_r$ for some sequence $\mathcal{G}$ constructed from $\mathcal{T}$. Note that each edge-maximal $K_5$-minor-free graph has a Wagner sequence by Theorem \ref{thm3.2}.

We actually can prove the following somewhat stronger result.
\begin{theorem}\label{thm3.5}
If $G$ is an edge-maximal $K_5$-minor-free graph, then $G$ has a 3-clique-coloring such that
 no triangle of $G$ is monochromatic.
\end{theorem}
\noindent{\bf Proof.}\,
If $|V(G)|\le 3$, the assertion is trivial. So let $|V(G)|\ge 4$. As
we have seen, $G$ has a Wagner sequence $\mathcal{T}$.
Let $\mathcal{T}=T_1, T_2, \ldots, T_r$ be a Wagner sequence of $G$.
We proceed by induction the {\em length} $r$ of $\mathcal{T}$.  When $r=1$, $G$ $(=T_1)$
is either the Wagner graph or a plane triangulation. If $G$ is the Wagner graph then the assertion is obvious, since the Wagner graph is 3-colorable.
If $G$ is a plane triangulation, then the assertion  follows directly from Lemma \ref{lem3.3}.
So assume that $r\ge 2$.

Note that $T_1, T_2, \ldots, T_{r-1}$ is  a  Wagner sequence
of the subgraph $G_{r-1}$ of $G$. By the induction hypothesis,  $G_{r-1}$  has a 3-clique-coloring such that
no triangle of $G_{r-1}$ is monochromatic. Let $\phi$
be such a 3-clique-coloring of $G_{r-1}$. It suffices to show that $\phi$
can be extended to a 3-clique-coloring of $G_r$ such that
no triangle of $G_{r}$ is monochromatic. Suppose that $T_r$ is pasted to $G_{r-1}$ along a
triangle $T$. Clearly the 3-clique-coloring $\phi$ of $G_{r-1}$ induces  a coloring of $\mathcal{H}(T)$, since $T$
is not monochromatic.
We can easily extend the coloring of $\mathcal{H}(T)$ to Int($A$) and to Ext($A$) (respectively) by applying Lemma \ref{lem3.4}. So the assertion follows.
Suppose that $T_r$ is pasted to $G_{r-1}$ along a $K_2$, say $A$. Then
$\phi$ induces a (not necessarily proper) coloring of $A$. If $A$
is  a maximal clique of $T_r$, then $T_r=V_8$. By using the 3-colorability of the Wagner graph $V_8$ and Lemma \ref{lem3.5}, we are done. Finally, if $A$
is not a maximal clique of $T_r$, we first extend  the  coloring
of $A$  to a 3-coloring of  $\mathcal{H}(T)$ where $T$ is the triangle of $T_r$ containing $A$,
and then extend the coloring of $\mathcal{H}(T)$  to Int($T$) and to Ext($T$) (respectively)
in $T_r$  by Lemma \ref{lem3.4}. Thus we obtain a 3-clique-coloring of $G_r$ such that no triangle of $G_{r}$ is monochromatic.
~$\Box$

\noindent{\bf Remark 2.}  By Theorem \ref{thm3.5}, we know that every edge-maximal $K_5$-minor-free graph $G$ is 3-clique-colorable.
Furthermore, we conjecture that this assertion is true for general $K_5$-minor-free graphs.


We now turn our attention to the claw-free graphs without $K_5$-minors.
Let $C_n+K_1$ be the graph obtained
from the disjoint union of $C_n$ and $K_1$ by joining the single vertex of $K_1$ to all the vertices of $C_n$.
The graph $C_n+K_1$ is also called a $n$-{\em wheel}, denoted by $W_n$,  and the vertex in $K_1$ is known as the
{\em hub} of $W_n$.

For claw-free graphs $G$ without $4$-cliques, we observe the following simple property
 of the graph $G$ by the Ramsey number $R(3,3)=6$,
its proof is similar to that of Lemma 8 in \cite{slk}.
\begin{lemma}\label{lem4.1}
If $G$ is a claw-free graph without $4$-cliques, then $\varDelta(G)\leq{5}$ and $G|N[v]$ is a 5-wheel $W_5$ for each vertex $v$ of degree $5$.
\end{lemma}

In \cite{slk} we proved that for a claw-free planar graph $G$, any 2-clique-coloring of $G-v$ can be extended to a 2-clique-coloring of $G$, where $v$  is a vertex of degree $5$ in $G$.  By Lemma \ref{lem4.1}, we can generalize
this result to \{claw, $K_5$-minor\}-free graphs. Its proof resembles that of Lemma 9 in \cite{slk},
and is omitted.

\begin{lemma}\label{lem4.2}
Let $G$ be a \{claw, $K_5$-minor\}-free graph without $4$-cliques and let $v$ be a vertex of degree $5$ in $G$. If $G-v$ is $2$-clique-colorable, then the same is true for $G$.\label{lem4.2}
\end{lemma}

\begin{lemma}\label{lem4.3}
Every \{claw, $K_5$-minor\}-free graph has maximum degree at most 6.
\end{lemma}
\noindent{\bf Proof.} Let $G$ be a \{claw, $K_5$-minor\}-free graph.
Suppose, to the contrary, that $\varDelta(G)\geq 7$. Let $v
\in V(G)$ such that $d_G(v)\geq 7$. Since $G$ contains no $K_5$-minor and claw, $G|N(v)$ contains no
$K_4$-minor, and so $\alpha(G|N(v))=2$. To obtain a contradiction,
we consider the graph $G|N(v)$.

Suppose that $G|N(v)$ contains a diamond $D$. Let $V(D)=\{u_1, u_2, u_3, u_4\}$ with $d_D(u_1)=d_D(u_2)=2$
and let $G_1=G|N(v)-V(D)$. Then $|V(G_1)|\ge 3$.
Since $\alpha(G|N(v))=2$, $G_1$ contains at most
two components and each vertex of $G_1$ is adjacent to $u_1$ or $u_2$.
Clearly, all vertices in each component of $G_1$ is adjacent to only one of $u_1$ and $u_2$, for otherwise
$G|N(v)$ would contain a $K_4$-minor.
If $G_1$ consists of precisely one component. Without loss of generality,
we may assume that all vertices of $G_1$ is adjacent to $u_1$. Thus
$u_2$ is not adjacent to any vertex of $G_1$. Observe that
$G|V(G_1)\cup\{u_1\}$ is not a complete subgraph in $G|N(v)$, since $|V(G_1)|\ge 3$ and  $G|N(v)$  contains no $K_4$-minor. Thus there exist veritces $u_5, u_6\in V(G_1)$ such that $u_5u_6\in E(G)$.  But then
$\{u_2, u_5, u_6\}$ is an independent set of $G|N(v)$, contradicting
the fact that $\alpha(G|N(v))=2$. If $G_1$ consists of precisely two component $O_1$ and $O_2$, and assume that $|V(O_1)|\ge |V(O_2)|$, then $|V(O_1)|\ge 2$. As we have observed above, all vertices of $O_i$ is adjacent
to exactly one of $u_1$ and $u_2$. Without loss of generality,
let us suppose that all vertices of $O_1$ is adjacent
to  $u_1$. So $u_2$ is not adjacent to any vertex of $O_1$. By $\alpha(G|N(v))=2$,
we see that $G|V(O_1)\cup\{u_1\}$  is complete, and thus $|V(O_1)|=2$,
 since $G|N(v)$ is $K_4$-minor-free. Hence $O_1=K_2$. On the other hand, we claim that all vertices of $O_2$ is adjacent
to  $u_2$. Indeed, if not, we take $c_i\in V(O_i)$ for $i=1,2$,
then $\{c_1,c_2, u_2\}$ is an independent set  of size 3 of
$G|N(v)$, a contradiction.  By the $K_4$-minor-freeness of $G|N(v)$,
it is easy to see that one of $u_3$ and $u_4$, say $u_3$, is not adjacent to any vertex of $O_2$.
$\alpha(G|N(v))=2$ implies that $u_3$ is adjacent to all vertices
of $O_1$. But then  $G|V(O_1)\cup\{u_1, u_3\}$ is a $K_4$-minor in $G|N(v)$, a contradiction.

Suppose that $G|N(v)$ contains no diamond. Let $I=\{u_1, u_2\}$ be a maximum independent set of $G|N(v)$, and let $N_1=N(v)-\{u_1, u_2\}$. As we have observed, each vertex of $N_1$ is adjacent to at least one of $u_1$ and $u_2$ by
$\alpha(G|N(v))=2$. If $G|N_1$ contains a triangle, then one of $u_1$ and $u_2$ is ajacent to
at least two vertices of this triangle. So $G|N(v)$ contains a diamond, contradicting our assumption.
Thus $G|N_1$ contains no triangle.
 By Ramsey number $R(3,3)=6$, we conclude that $|N_1|=5$ and
 $\alpha(G|N_1)=2$, and so $|N(v)|=7$. Clearly, $G|N_1$ is isomorphic to the cycle $C_5$.
 Note that
 either $u_1$ or $u_2$ is adjacent to
 at least three vertices on the cycle $G|N_1$. But then  $G|N(v)$ has a $K_4$-minor,
 a contradiction. ~~$\Box$

\begin{lemma}\label{lem4.4}
Let $G$ be a \{claw, $K_5$-minor\}-free graph with $\omega(G)=4$ and let $v$ be a vertex which lies in a $4$-clique
of $G$. If there exists a $2$-clique-coloring of $G-v$, then the same is true for $G$.
\end{lemma}

\noindent{\bf Proof.} The proof is by contradiction. Suppose that $G$ has no 2-clique-coloring. Let $\phi'$ be a 2-clique-coloring of $G-v$ with colors {\it red} and {\it green}. Then the extension of the coloring $\phi'$ of $G-v$ is impossible. Consequently, $G$ contains two maximal cliques $K$ and $L$ such that $V(K)\cap{V(L)={v}}$. Let $Q:=K-v$ and $R:=L-v$. Without loss of generality, we may assume that the vertices of $Q$ are red in $\phi'$, while those of $R$ are green. Thus we cannot color $v$ neither red nor green in any extension of $\phi'$.
Since $\phi'$ is a 2-clique-coloring of $G-v$, there exist two cliques (not necessarily maximal) $Q'$ and $R'$ in $G-v$ such that $Q'=Q+q_1$ and $R'=R+q_2$ with $q_1\not\in{V(Q)}$, $\phi'(q_1)=green$ and $q_2\not\in{V(R)}$, $\phi'(q_2)=red$, since otherwise $\phi'$ would not be a proper 2-clique-coloring of $G-v$.

Suppose that one of $K$ and $L$ is a $4$-clique of $G$. Without loss of generality,
let $K$ be a 4-clique of $G$. By the $K_5$-minor-freeness of $G$, clearly $q_1$ is not adjacent to $v$. Let us consider the graph $G-V(Q)$.
 If there is a path $P$ between $q_1$ and $v$ in $G-V(Q)$, then the vertices of $V(P)\cup V(Q)$  would contains a $K_5$-minor of $G$, contradicting our assumption.
 Hence $q_1$ and $v$ lie in the distinct components of $G-V(Q)$.
  Let $C$ be the component containing $v$. Let us now define a vertex coloring $\phi$ of $G$ as follows: we color $v$ green and change the colors of the vertices in $C-v$, and assign colors in $\phi'$ to all other vertices. We claim that $\phi$ is a 2-clique-coloring of $G$. Suppose not, let $M$ be a maximal clique of $G$ which is monochromatic in $\phi$. Then $M$ must contain at least one vertex, say $k$, of $Q$, and some vertices of $C$. This implies that $M$ is red, and thus $v\not\in{V(M)}$. Hence $M$ contains at least a vertex  $c$  that
 is not adjacent to $v$. One can easily see that $c\in{V(C)}$.
 Note that there is no  path between $q_1$ and $v$ in $G-V(Q)$,
 so $c$ is not adjacent to $q_1$. But then we find  a claw induced by $\{k,q_1,v,c\}$ centered at $k$, a contradiction. Thus, we may assume that neither $K$ nor $L$ is a $4$-clique of $G$.
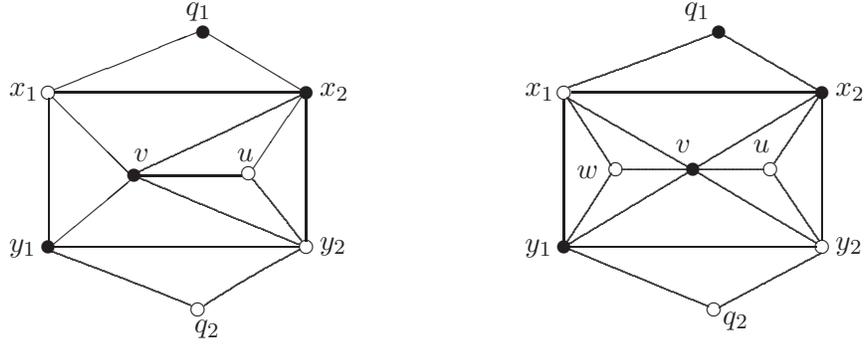
\begin{figure}[ht]\label{map1}
\setlength{\unitlength}{.045in}
\begin{picture}(55,35)
\put(30,-6){\begin{picture}(60,35)

\put(-4,28){\circle{1.5}}   \put(-4,10){\circle*{1.5}}   \put(26,28){\circle*{1.5}}   \put(26,10){\circle{1.5}}   \put(6,18.3){\circle*{1.5}}
\put(19.3,18.6){\circle{1.5}}   \put(14,35){\circle*{1.5}}   \put(13.4,2.7){\circle{1.5}}

\put(-3.3,28){\line(1,0){29.7}}   \put(-3.9,27.3){\line(0,-1){17.7}}
\put(26,28){\line(0,-1){17.3}}  \put(-4,10){\line(1,0){29.2}} \put(-3.5,27.5){\line(1,-1){9.5}}
\qbezier(25.5,27.7)(16,23)(6,18.5)   \put(-4,9.8){\line(6,5){10}}
\qbezier(6,18.3)(15,14.4)(25.2,10.3)
\put(26,28.3){\line(-2,-3){6.1}}

\qbezier(25.5,10.7)(23,14)(19.8,17.9)   \put(-3.3,28.3){\line(5,2){17.5}}
\put(26,28){\line(-5,3){11.5}}   \qbezier(-4.2,9.7)(4.8,6.4)(12.6,3.2)
\qbezier(14.2,3.2)(19,6.4)(25.3,9.7)
\put(18.5,18.3){\line(-1,0){13}}

\put(-8.5,27.6){$x_1$}  \put(27.6,27.6){$x_2$}  \put(12,37){$q_1$}
\put(-8.5,9.5){$y_1$}   \put(27.6,9.7){$y_2$}  \put(6,20){$v$}  \put(18,20){$u$}  \put(13,0){$q_2$}

\put(56,28){\circle{1.5}}   \put(56,10){\circle*{1.5}}   \put(86,28){\circle*{1.5}}   \put(86,10){\circle{1.5}}
\put(71,19){\circle*{1.5}}
\put(80,19){\circle{1.5}}   \put(62,19){\circle{1.5}}   \put(74,35){\circle*{1.5}}   \put(73.4,2.7){\circle{1.5}}

\put(56.9,28){\line(1,0){29.7}}   \put(56,27.13){\line(0,-1){17.7}}
\put(86,28.3){\line(0,-1){17.5}}  \put(56,10){\line(1,0){29.3}} \qbezier(56.5,27.3)(63.5,23.25)(71,19)
\qbezier(71,19)(78.5,23.5)(86,28)   \qbezier(56,10)(63.5,14.5)(71,19)    \qbezier(71,19)(78.5,14.5)(85.4,10.44)   \qbezier(80.5,19.8)(83,23.5)(86,28)  \qbezier(80.5,18.4)(83.1,14.55)(85.7,10.7)   \qbezier(56.6,28.6)(65,31.5)(74,35)
\qbezier(74,35)(80,31.5)(86,28)   \qbezier(56,10)(64.7,6.35)(72.8,3.1)   \qbezier(74.1,3)(80,6.125)(85.9,9.25)
\qbezier(71,19)(75.5,19)(79.2,19)  \qbezier(56.4,27.4)(59,23.5)(61.5,19.5)   \qbezier(71,19)(66.5,19)(62.8,19)
\qbezier(56,10)(59,14.5)(61.5,18.5)

\put(51.5,27.6){$x_1$}  \put(87.6,27.6){$x_2$}  \put(70,37){$q_1$}
\put(51.5,9.5){$y_1$}   \put(87.6,9.7){$y_2$}  \put(69,21){$v$}  \put(78,21){$u$}  \put(74.5,1){$q_2$}
\put(57.5,18){$w$}

\end{picture}}
\end{picture}
\\

\caption{The graphs $F_1$ and $F_1'$.}
\end{figure}

To complete the proof, we have the following claim.

{\bf Claim 1} Both $K$ and $L$ are $3$-cliques  of $G$, and the induced subgraph $G|N[v]\cup\{q_1, q_2\}$ is isomorphic to  $F_1$ or $F_1'$ (see Figure 2).

We first show that both $K$ and $L$ are $3$-cliques  of $G$.
Indeed, if not, without loss of generality, we may assume that $K$ be a 2-clique of $G$ and let $V(K)=\{v, x\}$. Obviously,
$x$ is not adjacent to any neighbor of $v$  by the maximality of $K$.
According to our assumption, the vertex $v$ lies in a $4$-clique, say $W$, of $G$.
So $x$ is not adjacent to any vertex of $L-v$ ($=R$) and $W-v$.
On the other hand, since $L$ and $W$ are two distinct cliques of $G$, there exists
vertices $y\in V(L)-V(W)$ and $u\in V(W)-V(L)$ such that $yu\not\in E(G)$. This
implies that $\{v, x, u, y\}$ induces a claw centered at $v$, a contradiction. Therefore,
$K$ and $L$ are $3$-cliques  of $G$.
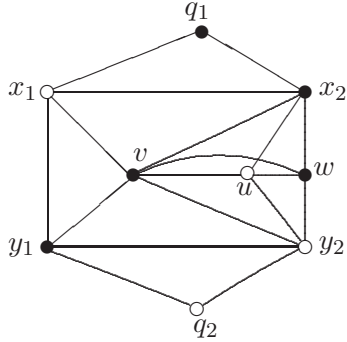
\begin{figure}[ht]\label{map1}
\setlength{\unitlength}{.045in}
\begin{picture}(55,40)
\put(55,-0){\begin{picture}(60,40)

\put(-4,28){\circle{1.5}}   \put(-4,10){\circle*{1.5}}   \put(26,28){\circle*{1.5}}   \put(26,10){\circle{1.5}}   \put(6,18.3){\circle*{1.5}}
\put(19.3,18.6){\circle{1.5}}   \put(14,35){\circle*{1.5}}   \put(13.4,2.7){\circle{1.5}}

\put(-3.3,28){\line(1,0){29.7}}   \put(-3.9,27.3){\line(0,-1){17.7}}
  \put(-4,10){\line(1,0){29.2}} \put(-3.5,27.5){\line(1,-1){9.5}}
\qbezier(25.5,27.7)(16,23)(6,18.5)   \put(-4,9.8){\line(6,5){10}}
\qbezier(6,18.3)(15,14.4)(25.2,10.3)
\put(26,28.3){\line(-2,-3){6.1}}

\qbezier(25.5,10.7)(23,14)(19.8,17.9)   \put(-3.3,28.3){\line(5,2){17.5}}
\put(26,28){\line(-5,3){11.5}}   \qbezier(-4.2,9.7)(4.8,6.4)(12.6,3.2)
\qbezier(14.2,3.2)(19,6.4)(25.3,9.7)
\put(18.5,18.3){\line(-1,0){13}}

\qbezier(6,18.3)(16,23)(26,18.3)\put(26,18.3){\circle*{1.5}} \put(27,18.3){$w$}
\qbezier(26,18.3)(26,23)(26,28)\qbezier(26,18.3)(26,13)(26,10.7)
\qbezier(26,18.3)(23,18.3)(20,18.3)

\put(-8.5,27.6){$x_1$}  \put(27.6,27.6){$x_2$}  \put(12,37){$q_1$}
\put(-8.5,9.5){$y_1$}   \put(27.6,9.7){$y_2$}  \put(6,20){$v$}  \put(18,16){$u$}  \put(13,0){$q_2$}

\end{picture}}
\end{picture}
\caption{The graph $H_1$.}
\end{figure}

Let $V(K)=\{v,x_1,x_2\}$, $V(L)=\{v,y_1,y_2\}$ and let $W$ be a 4-clique that contains the vertex $v$. Note that
 $v\in V(W)\cap V(K)\cap V(L)$, so
 $1\le|V(W)\cap(V(K)\cup V(L))|\le 3$, for otherwise $V(W)\supseteq V(K)$ or $V(W)\supseteq V(L)$,
 contradicting the fact that $K$ and $L$ are maximal cliques of $G$. Thus
 $V(W)-(V(K)\cup V(L))\neq\emptyset$. Let $u\in V(W)-(V(K)\cup V(L))$. Clearly,
 $K$ (respectively $L$) contain at least a vertex that is not adjacent to $u$. Without loss of generality,
 let $x_1\in V(K), y_1\in V(L)$ such that $x_1u\not\in E(G)$ and $y_1u\not\in E(G)$. Since
 $\{v, x_1, y_1, u\}$ does not induce a claw centered at $v$, it immediately follows that
 $x_1y_1\in E(G)$.
This implies that $x_1y_2\not\in E(G)$ and $x_2y_1\not\in E(G)$ by the maximality of $K$ and $L$.
Since $\{v, x_1, y_2, u\}$ and $\{v, x_2, y_1, u\}$ can not induce a claw centered at $v$, we have $ux_2, uy_2\in E(G)$. Now we show that $G|N[v]\cup\{q_1, q_2\}$ is isomorphic to either $F_1$
or $F_1'$. If $N(v)=\{x_1, x_2, y_1, y_2, u\}$, then $x_2y_2\in E(G)$ since $v$ lies in
a 4-clique of $G$. Hence $v$ lies in the 4-clique induced by $\{v, u, x_2, y_2\}$.
Furthermore, we claim that $uq_i\not\in E(G)$ for $i=1,2$.
If not, let $uq_1\in E(G)$, then $G|N[v]\cup \{q_1\}$ would contain a $K_5$-minor.
Therefore, $G|N[v]\cup\{q_1, q_2\}$ is isomorphic to the graph $F_1$.
If $N(v)\neq\{x_1, x_2, y_1, y_2, u\}$,
we have $d_G(v)=\varDelta(G)=6$ by Lemma \ref{lem4.3}. Let
 $N(v)=\{x_1, x_2, y_1, y_2, u, w\}$.
 Since $\{v, x_1, y_2, w\}$ does not induce a claw centered at $v$,
  $wx_1\in E(G)$ or  $wy_2\in E(G)$.
  If $wy_2\in E(G)$, then $wy_1\not\in E(G)$ by the maximality of $L$. So $wx_2\in E(G)$ (see the graph $H_1$ in Figure 3), for otherwise $\{v, x_2, y_1, w\}$ would induce a claw centered at $v$. This implies that
 $wx_1\not\in E(G)$ by the maximality of $K$. As we have seen, $u$ is not adjacent to $x_1$, it follows that
$wu\in E(G)$, since otherwise $\{v, x_1, u, w\}$ induces a claw centered at $v$. But now
$G[N[v]]$ contains a $K_5$-minor, a contradiction.
If $wx_1\in E(G)$, then $wx_2 \not\in E(G)$ by the maximality of $K$. To avoid a claw  induced by
$\{v, x_2, y_1, w\}$ centered at $v$, we have $wy_1\in E(G)$. Thus $wy_2 \not\in E(G)$ by the maximality of $L$. By
the $K_5$-minor-freeness of $G$, it is easy to see that $wu\not\in E(G)$.
Note that $\{v,w, x_2, y_2\}$ can not induce a claw centered at $v$, so $x_2y_2\in E(G)$.
Finally, one easily see that $uq_i\not\in E(G)$ and $wq_i\not\in E(G)$ by the $K_5$-minor-freeness of $G$.
Consequently, $G|N[v]\cup\{q_1, q_2\}$ is isomorphic to the graph $F_1'$.
~~$\Box$

For convenience, let us denote by $W$ and $W'$ (if exists) the 4-cliques $G|\{x_2,y_2, u,v\}$ and $G|\{x_1,y_1, v,w\}$, that contain the vertex $v$, in $F_1$ or $F_1'$.

 As we saw earlier, $\phi'(x_1)=\phi'(x_2)=red$, $\phi'(y_1)=\phi'(y_2)=green$, and $\phi'(q_1)=green$ and $\phi'(q_2)=red$. We give a 2-clique-coloring $\phi$ of $G$ as follows: we exchange the colors of $x_2$ and $y_2$, and assign red or green to $v$, and let $\phi(x)=\phi'(x)$ for all the vertices $x\in{V(G)-\{v,x_2,y_2\}}$. We claim that $\phi$ is a 2-clique-coloring of $G$. Suppose not, let $M$ be a monochromatic maximal clique of $G$ in $\phi$. Then $M$ must contain exactly one vertex of $\{x_2,y_2\}$, and at least one vertex, say $k$, of $V(G)-(N[v]\cup \{q_1,q_2\})$. By the symmetry between $x_2$ and $y_2$ in $F_1$ or $F_1'$, we may assume
  that $x_2$ is in $M$. This implies that $M$ is green. By Lemma \ref{lem4.3},
 $d_G(x_2)\le \varDelta(G)\le 6$. If $|V(M)|\ge 3$, then
 all the vertices in $V(M)-\{k\}$ are in $N[v]\cup \{q_1,q_2\}$; otherwise we have
 $$d_G(x_2)=|N(x_2)|
 \ge|\{x_1,y_2, u, v, q_1\}|+|V(M)-(N[v]\cup \{q_1,q_2\})|\ge 7,$$
which is a contradiction. We claim that $|V(M)|\le 3$.  Indeed,
if $|V(M)|=4$, by Claim 1, $G|N[v]\cup \{q_1,q_2\}\cong F_1$ or $F_1'$, then $k$ is not adjacent to
$v$. Since $M$ is a  maximal green 4-clique of $G$, we have $V(M)-\{k\}=\{x_2, u, q_1\}$.
But now we can find four vertex-disjoint paths linking
$x_1$ to all vertices in the 4-clique $W$. Thus $G$ contains a $K_5$-minor, contradicting our assumption.
Consequently, $M$ is  either a green 2-clique or a green 3-clique of $G$. We consider the following two cases.

{\em Case 1}: $M$ is a $2$-clique of $G$, that is, $M$ is the maximal 2-clique induced by $\{x_2,k\}$.  Obviously, we have $x_1k\not\in{E(G)}$, $y_2k\not\in{E(G)}$ by the maximality of $M$. By Claim 1, we know that $x_1y_2\not\in{G}$. This implies that  $\{x_1,y_2,k,x_2\}$ induces a claw centered at $x_2$, a contradiction.

\begin{figure}[ht]\label{map1}
\setlength{\unitlength}{.043in}
\begin{picture}(500,40)
\put(60,-6){\begin{picture}(500, 40)

\put(-40,28){\circle{1.5}}   \put(-40,10){\circle*{1.5}}   \put(-10,28){\circle*{1.5}}   \put(-10,10){\circle{1.5}}   \put(-30,18.3){\circle*{1.5}}
\put(-16.7,18.6){\circle*{1.5}}   \put(-22,35){\circle*{1.5}}   \put(-22.6,2.7){\circle{1.5}}
\put(5,33){\circle*{1.5}}
\put(-36,18.3){\circle*{1.5}}

\put(-39.3,28){\line(1,0){29.7}}   \put(-40,27.3){\line(0,-1){17.7}}
\put(-10,28){\line(0,-1){17.6}}  \put(-40,10){\line(1,0){29.35}} \put(-39.6,27.6){\line(1,-1){9.8}}
\put(-10,28){\line(-2,-1){20}}   \put(-40,10){\line(6,5){10}}    \put(-10.5,10.25){\line(-5,2){19.9}}   \put(-10,28){\line(-2,-3){6}}  \put(-10.5,10.3){\line(-3,4){5.8}}   \put(-39.6,28.2){\line(5,2){17.5}}
\put(-10,28){\line(-5,3){11.5}}   \put(-40,9.8){\line(5,-2){17}}    \put(-10.45,9.8){\line(-5,-3){11.5}}
\put(-17.2,18.3){\line(-1,0){13}}
\put(-10,28){\line(3,1){15}}
\qbezier(4.5,32.5)(-5,25.5)(-16,18.4)

\qbezier[20](-36.5,18.3)(-33,18.3)(-31,18.3)
\qbezier[30](-36,18.5)(-37.8,23)(-39.7,27.3)
\qbezier[30](-36,18.3)(-37.8,14)(-39.9,10.2)

\put(-44.5,27.6){$x_1$}  \put(-10,31){$x_2$}  \put(-22,37){$q_1$}
\put(-44.5,9.5){$y_1$}   \put(-8.4,9.7){$y_2$}  \put(-30,21){$v$}  \put(-19,21){$u$}  \put(-22,0){$q_2$}
   \put(7,31){$k$} \put(-39.5,18){$w$}

\put(25,28){\circle{1.5}}   \put(25,10){\circle*{1.5}}   \put(55,28){\circle*{1.5}}   \put(55,10){\circle{1.5}}   \put(35,18.3){\circle*{1.5}}
\put(48.3,18.6){\circle{1.5}}   \put(43,35){\circle*{1.5}}   \put(42.4,2.7){\circle{1.5}}
\put(68,34.7){\circle*{1.5}}
\put(29,18.3){\circle*{1.5}}

\put(25.6,28){\line(1,0){29.7}}   \put(25,27.4){\line(0,-1){17.7}}
\put(55,28){\line(0,-1){17.5}}  \put(25,10){\line(1,0){29.3}} \put(25.3,27.5){\line(1,-1){9.55}}
\put(55,28){\line(-2,-1){20}}   \put(25,10){\line(6,5){10}}    \put(54.5,10.35){\line(-5,2){20}}   \put(55,28){\line(-2,-3){6}}  \put(54.8,10.5){\line(-3,4){5.9}}   \put(25.5,28.3){\line(5,2){17.5}}
\put(55,28){\line(-5,3){11.4}}   \put(25,10){\line(5,-2){17.1}}    \put(54.4,9.8){\line(-5,-3){11.5}}
\put(47.9,18.3){\line(-1,0){13}}
\put(55,28){\line(2,1){13}}  \put(68,34.9){\line(-1,0){24.5}}
\qbezier[20](29,18.3)(33,18.3)(34.5,18.3)\qbezier[30](29,18.3)(27,22.5)(25.2,27.2)
\qbezier[30](29,18.3)(27,14)(25,10.2)

\put(20.5,27.6){$x_1$}  \put(56,26.5){$x_2$}  \put(43,37){$q_1$}
\put(20.5,9.5){$y_1$}   \put(56.6,9.7){$y_2$}  \put(35,21){$v$}  \put(46,21){$u$}  \put(43,0){$q_2$}
\put(67,31){$k$}\put(25,18){$w$}

\put(41,-5){$F_3$}   \put(-24,-5){$F_2$}

\end{picture}}
\end{picture}
\\
\\

\caption{The graphs $F_2$ and $F_3$.}
\end{figure}
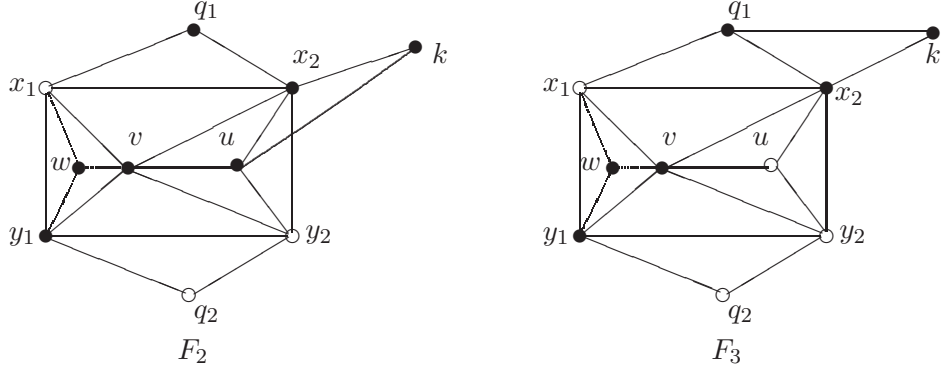

{\em Case 2}: $M$ is a 3-clique of $G$, and let $V(M)=\{x_2, k,l\}$.

As we have seen,  $l\in N[v]\cup \{q_1,q_2\}$, and $k$ is not adjacent to
$v$ by Claim 1. Since $M$ is a  maximal green 3-clique of $G$, we have $l=u$ or $q_1$.

 If $l=u$, that is, $ku\in E(G)$, then $\phi(u)=green$ (see the graph $F_2$ in Figure 4). This implies that $y_2k\not\in E(G)$ by the maximality of $M$.  Thus $x_1k\in E(G)$ for avoiding the claw $G|\{x_1, x_2, y_2, k\}$ at $x_2$. Hence we can now find four vertex-disjoint paths linking
$x_1$ to all vertices in the 4-clique $W$, and so a $K_5$-minor occurs in $G$, a contradiction.

If $l=q_1$, that is, $q_1k\in E(G)$, then $V(M)=\{x_2, k, q_1\}$ (see the graph $F_3$ in Figure 4). By the maximality of $M$,
we have $x_1k\not\in E(G)$. By Claim 1, we know that $x_1y_2\not\in E(G)$ and $x_1u\not\in E(G)$.
Note that $\{x_1,x_2, y_2, k\}$ and $\{x_1,x_2, u, k\}$ can not induce
claws at $x_2$, so $y_2k, uk\in E(G)$, so we find four vertex-disjoint paths linking
$k$ to all vertices in the 4-clique $W$. But now this produces a $K_5$-minor  in $G$, a contradiction.
~~$\Box$

\begin{lemma}[Bacs\'o and Tuza \cite{bz}]\label{lem4.5}
 Every connected claw-free graph of maximum degree at most four, other than an odd
cycle, is 2-clique-colorable. Moreover, a 2-clique-coloring can be found in polynomial time.
\end{lemma}

Finally, we have the following result.

\begin{theorem}\label{thm4.1}
Every \{claw, $K_5$-minor\}-free graph $G$ of order $n$, different from an odd cycle, is $2$-clique-colorable.
\end{theorem}

\noindent{\bf Proof.} We proceed by induction on $n$. For $n\leq{4}$, clearly the assertion holds. Now let $n>{4}$, and assume that the assertion holds for smaller values than $n$. If $\varDelta(G)\le 4$, by Lemma \ref{lem4.5}, $G$ is $2$-clique-colorable.
So we may assume that $\varDelta(G)\ge 5$.

Suppose that $G$ has no $4$-cliques. Then $\varDelta(G)\leq{5}$ by Lemma \ref{lem4.1}, and so $\varDelta(G)=5$.
 Let $v$ be a vertex of degree $5$ in $G$. By Lemma \ref{thm4.1}, we have $G|N[v]$ is a $5$-wheel. If $G=G|N[v]$,
   then $G$ can easily be clique-colored in two colors. If $G\neq G|N[v]$, then clearly $G-v$ is still a \{claw, $K_5$-minor\}-free graph, not an odd cycle. Therefore, by the induction hypothesis, $G-v$ is 2-clique-colorable.
It follows from Lemma \ref{lem4.2} that $G$ is 2-clique-colorable.

Suppose that $G$ has $4$-cliques. Let $v$ is a vertex of $G$ that lies in a $4$-clique. Obviously $G-v$ is still a \{claw, $K_5$-minor\}-free graph, not an odd cycle. By the induction hypothesis, $G-v$ is $2$-clique-colorable. So $G$ is also 2-clique-colorable by Lemma \ref{lem4.4}.
~$\Box$

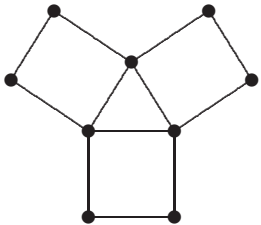
\begin{figure}[ht]\label{map1}
\setlength{\unitlength}{.045in}
\begin{picture}(60,26)
\put(50,-3){\begin{picture}(60,26)

\multiput(10,10)(10,0){2}{\circle*{1.5}}
\multiput(10,0)(10,0){2}{\circle*{1.5}}
\put(15,18){\circle*{1.5}}
\multiput(1,16)(28,0){2}{\circle*{1.5}}
\multiput(6,24)(18,0){2}{\circle*{1.5}}

\multiput(10,0)(10,0){2}{\line(0,1){10}}
\multiput(10,0)(0,10){2}{\line(1,0){10}}

\qbezier(15,18)(12.5,14)(10,10)\qbezier(15,18)(17.5,14)(20,10)
\qbezier(15,18)(10.5,21)(6,24)\qbezier(15,18)(19.5,21)(24,24)

\qbezier(1,16)(5.5,13)(10,10)\qbezier(1,16)(3.5,20)(6,24)

\qbezier(29,16)(24.5,13)(20,10)\qbezier(29,16)(26.5,20)(24,24)
\end{picture}}
\end{picture}
\\

\caption{Example of a $K_5$-minor-free graph containing claws that is not 2-clique-colorable.}
\end{figure}

\noindent{\bf Remark 3.} The condition \{claw, $K_5$-minor\}-free in Theorem \ref{thm4.1}
cannot be dropped. For example, the graph shown in Figure 4 contains a claw and its clique-chromatic number
is 3. The line graph $L(K_6)$ of $K_6$ contains the complete graph $K_5$, and is not 2-clique-colorable by Ramsey number $R(3, 3)=6$.

Note that if $\phi$ is a $2$-clique-coloring of a graph, then
$\phi^{-1}(r)$ and $\phi^{-1}(g)$ are clique-transversal sets of $G$. By Theorem \ref{thm4.1}, we
immediately obtain an  upper bound on the clique-transversal number for \{claw, $K_5$-minor\}-free graphs.

\begin{corollary}\label{co}
Every \{claw, $K_5$-minor\}-free graph, different from an odd cycle, has the clique-transversal number
 bounded above by half of its order.
\end{corollary}

 Finally, we present a polynomial-time algorithm to find a 2-clique-coloring of \{claw, $K_5$-minor\}-free graphs. In \cite{bz} Bacs\'o and Tuza proposed the polynomial-time algorithm CLQCOL for 2-clique-coloring problem on claw-free graphs of maximum degree at most four, other than an odd hole.

Clearly, if $G$ is a \{claw, $K_5$-minor\}-free graph,  not an odd cycle,  then so is the graph $G-v$, and
$G-v$ has fewer vertices. Based on the algorithm CLQCOL and Lemmas 3.4-3.7, we provide the following algorithm to find a 2-clique-coloring on \{claw, $K_5$-minor\}-free graphs different from odd cycles.


 \noindent{\bf Algorithm} $\mathcal{A}$. 2-Clique-coloring of \{claw, $K_5$-minor\}-free graphs.

{\em Input:} \{Claw, $K_5$-minor\}-free graph $G$, not an odd cycle.

{\em Output:} 2-Clique-coloring $\phi$: $V(G)\rightarrow\{r,g\}$.

{\it Step 1:} If $\Delta(G)\le{4}$, then CLQCOL($G$)(see, Ref. \cite{bz}), stop the algorithm.
If not, turn to Step 2.

{\em Step 2:} If $\Delta(G)={5}$ and there is no 4-clique in $G$, turn to Step 3.
If not, turn to Step 4.

{\em Step 3:} If $G$ is a 5-wheel, give a 2-clique-coloring directly. If not, find a vertex $v$ with degree 5. Then $\mathcal{A}(G-v)$. Extend the 2-clique-coloring of $G-v$ to a 2-clique-coloring of $G$. Stop the algorithm.

{\em Step 4:} Find a vertex $v$ of degree $\le{5}$ in $G$ such that $v$ lies in a 4-clique. Then $\mathcal{A}(G-v)$. Extend the 2-clique-coloring of $G-v$ to a 2-clique-coloring of $G$. Stop the algorithm.

 By Lemmas \ref{lem4.2}, \ref{lem4.4} and \ref{lem4.5}, the validity
of Algorithm $\mathcal{A}$ can be easily verified. It is easy to check that the loop of every step is performed at
most $n$ times. In Step 1, the running time of Algorithm CLQCOL is $O(n^2)$ and is carried out only once.
In Steps 2, 3 and 4, the running time is at most $O(n)$. Thus the overall time is at most $O(n^2)$. Consequently,
we obtain the following result.

\begin{theorem}
Algorithm $\mathcal{A}$ is a polynomial-time algorithm in $O(n^2)$ for the 2-clique-coloring of \{claw, $K_5$-minor\}-free graphs, different from odd cycles.
\end{theorem}

\section*{Acknowledgment}
The research was  supported in part by grant 11171207 of the
National Nature Science Foundation of China.

\end{document}